\documentclass[11pt,a4paper]{article}
\usepackage[english]{babel}
\usepackage{amsmath, amsthm, amsfonts, amssymb, amscd}
\usepackage[utf8]{inputenc}
\usepackage[pdftex]{graphicx}
\usepackage[matrix,arrow]{xy}
\usepackage{enumerate}
\usepackage{xcolor}
\usepackage{comment}
\DeclareGraphicsExtensions{.bmp,.png,.pdf,.jpg}

\newtheorem{teo}{Theorem}
\newtheorem{lema}[teo]{Lemma}
\newtheorem{cor}[teo]{Corollary}

\newtheorem{rem}{Remark}

\newcommand{\pf}{\noindent{\bf Proof\ \ }}
\newcommand{\cqd}{{\hfill $\rule{2mm}{2mm}$}\vspace{3mm}}
\newcommand{\Z}{\mathbb Z}
\newcommand{\N}{\mathbb N}

\newcommand{\I}{\mathcal{I}}
\newcommand{\J}{\mathcal{J}}

\newcommand{\f}{\mathfrak f}
\newcommand{\Iff}{\Longleftrightarrow}
\usepackage{float}

\begin{document}
\title{On value sets of fractional ideals}

\author{E. M. N. de Guzm\'an\footnote{Supported by a fellowship from CAPES}  \and A. Hefez\footnote{Partially supported by the CNPq Grant 307873/2016-1} }

\date{ }
\maketitle

\noindent{\bf Abstract}
The aim of this work is to study duality of fractional ideals with respect to a fixed ideal and to investigate the relationship between value sets of  pairs of dual ideals in admissible rings, a class of rings that contains the local rings of algebraic curves at singular points. We characterize canonical ideals by means of a symmetry relation between lengths of certain quotients of associated ideals to a pair of dual ideals. In particular, we extend the symmetry among absolute and relative maximals in the sets of values of pairs of dual fractional ideals to other kinds of maximal points. Our results generalize and complement previous ones by other authors.
\bigskip

\medskip

\section{Introduction}

Value sets or semigroups of rings of irreducible plane curves germs, called plane branches, were studied by Zariski in \cite{Za1} and their importance is due to the fact that they constitute, over $\mathbb C$, a complete set of discrete invariants for their topological classification. From the work of Ap\'ery \cite{Ap}, it follows that this semigroup in the set $\N$ of natural numbers is in some sense symmetric. Many years later, Kunz, in \cite{Ku}, motivated by a question asked by Zariski, showed that a one dimensional noetherian domain, with some additional technical conditions, is Gorenstein if and only if its semigroup of values is symmetric.

For a germ of a singular plane curve with several branches over $\mathbb C$, Waldi in \cite{W}, based on the work \cite{Za1} of Zariski, showed that also in this case the topological type of the germ is 
characterized by the semigroup of values of the local ring of the curve, this time, a semigroup of $\N^r$, where $r$ is the number of branches of the curve.  Although not finitely generated, this semigroup was shown by Garcia in \cite{Ga}, for $r=2$, to be determined in a combinatorial way by a finite set of points that he called \emph{maximal points}. Garcia also showed that these maximal points of the semigroup of a plane curve have a certain symmetry. These results were generalized later, for any value of $r$, by Delgado in \cite{D87}, where two kinds of maximal points were emphasized:  the relative and absolute maximals, showing that the relative maximals determine the semigroup of values in an inductive and combinatorial way and that the relative and absolute maximals determine each other, extending Garcia's symmetry. A short time later, Delgado, in \cite{D88}, generalizing the work of Kunz, introduced a concept of symmetry for value semigroups in $\N^r$ and showed that this symmetry is equivalent to the Gorensteiness of the ring of the curve.

In \cite{Da}, D'Anna, generalizing the work of J\"ager, in \cite{ja}, and of Campillo, Delgado and Kiyek, in \cite{CDK}, extended the properties of value semigroups for some class of one dimensional noetherian rings to value sets of their regular fractional ideals and characterized a normalized canonical ideals of a given ring in terms of a precisely described value set obtained from the value semigroup of the ring.

Also, recently, Pol, in the work \cite{Pol17}, extended Delgado's result by showing that the Gorensteiness of the ring of a singularity is equivalent to some symmetry relation among sets of values of any dual pair of regular fractional ideal, duality taken with respect to the ring itself, and also showed that in this case one has a pairing between absolute and relative maximal points of the value set of an ideal and that of its dual. In the work \cite{KST} (see also \cite{Pol18}), the authors show that this symmetry relation among value sets of dual pairs of ideals holds without any extra assumption on the ring, if one takes duality with respect to a canonical ideal.

In this paper, we generalize the work \cite{CDK} that characterizes the Gorensteiness of an admissible ring in terms of lengths of certain quotients of complementary ideals by establishing similar conditions, valid without the Gorenstein assumption, for pairs of dual regular fractional ideals with respect to any fixed fractional ideal. This will allow us to unify, generalize and complement previous results in \cite{Pol17}, \cite{Pol18} and \cite{KST} and get new symmetry relations among other types of maximal points, other than absolute and relative maximal points, in value sets of pairs of dual fractional ideals. 

\section{Admissible rings, fractional ideals and value sets}

Following \cite[Definition 3.5]{KST}, a noetherian, Cohen-Macauley, one dimensional local ring $(R,\mathfrak M)$ is called \emph{admissible}, if it is analytically reduced, residually rational and  $\# R/\mathfrak M \geq r$, where $r$ is the number of valuation rings over $R$ of the total ring of fractions $Q$ of $R$. In this context, all valuation rings are discrete valuation rings (cf. \cite[Theorem 3.1]{KST}).  

This class of rings, without any special given name, was previously considered in \cite{CDK} and in \cite{Da}, generalizing the important family of local coordinate rings of reduced curves at a singular point.

Throughout this work we will assume that $R$ is an admissible ring. We denote by $\Z$ the set of integer and by $I$ the set 
$\{1,\ldots,r\}$. If $v_1,\ldots,v_r$ are the valuations associated to the discrete valuation rings of $Q$ over $R$, then we have a value map
$v\colon Q^{reg} \to \Z^r$, where $Q^{reg}$ is the set of regular elements of $Q$, defined by $h\mapsto v(h)=(v_1(h),\ldots,v_r(h))$  (cf. \cite[Definition 3.2]{KST}). We will consider on $\Z^r$ the natural partial order $\leq$ induced by the order of $\Z$. 

An $R$-submodule $\mathcal I$ of $Q$ will be called a \emph{fractional ideal} if there is a regular element $d$ in $R$ such that $d\,\mathcal I \subset R$. The ideal $\I$ will be said a \emph{regular fractional ideal}, if it contains a regular element of $Q$.

Examples of regular fractional ideals of $R$ are $R$ itself, the integral closure $\widetilde{R}$ of $R$ in $Q$, any ideal of $R$ or of $\widetilde{R}$ that contains a regular element and the ideals of the form $$\J \colon \I=\{x\in Q; \ x\I \subset \J\},$$ where $\mathcal I$ and $\mathcal J$ are regular fractional ideals. In particular, the conductor $\mathcal C=R\colon \widetilde{R}$, which is the largest common ideal of $R$ and $\widetilde{R}$, is a regular fractional ideal. Notice that $\J\colon R=R$ for all fractional ideal $\J$.

In the class of admissible rings one has that there is a natural isomorphism $\mathcal J \colon \mathcal I \simeq Hom_R(\mathcal I,\mathcal J)$ (cf. \cite[Lemma 2.4]{KST}), for any fractional ideals $\I$ and $\J$. It is always true that 
$\I \subseteq \J\colon(\J\colon \I)$. The fractional ideal $\J$ is called a \emph{canonical ideal} if the last inclusion is an equality for every fractional ideal $\I$. In our context, canonical ideals exist (cf. \cite{Da} or \cite{KST}); two canonical ideals differ up to a multiplication by a unit in $Q$ and a multiple by such a unit of a canonical ideal is also a canonical ideal. The ring $R$ will be called \emph{Gorenstein} if $R$ itself is a canonical ideal.

We define the \emph{value set} of a regular fractional ideal $\I$ of $R$ as being 
$$E(\I)=v(\mathcal I^{reg})\subset \Z^r.$$

The value set $E(R)$ of $R$ is a subsemigroup of $\N^r$, called the \emph{semigroup of values} of $R$. The value set $E(\mathcal I)$ of a fractional ideal $\mathcal I$ is not necessarily closed under addition, but it is such that $E(R)+E(\mathcal I)\subset E(\mathcal I)$. For this reason $E(\I)$ is called a \emph{semigroup ideal} of the semigroup $E(R)$. More generally, one has
\begin{equation}\label{$E+E^*$}
E(\I)+E(\J\colon \I)\subset E(\J).
\end{equation}

A value set $E$ of a regular fractional ideal has the following fundamental properties (cf. \cite[Proposition 3.9]{KST}):

\noindent $E_0$: \ 
There are $\alpha\in \Z^r$ and $\beta \in \N^r$ such that $\beta+\N^r\subset E \subset \alpha+\N^r$;

\noindent $E_1$: \	If $\alpha=(\alpha_1,\ldots,\alpha_r)$ and $\beta=(\beta_1,\ldots,\beta_r)$ belong to $E$, then 
	$$\min(\alpha,\beta)=(\min(\alpha_1,\beta_1),\ldots,\min(\alpha_r,\beta_r))\in E;$$
\noindent $E_3$: \ If $\alpha=(\alpha_1,\ldots,\alpha_r), \beta=(\beta_1,\ldots,\beta_r)$ belong to $E$, $\alpha\neq\beta$ and $\alpha_i=\beta_i$ for some $i\in\{1,\ldots,r\}$, then there exists $\gamma\in E$ such that $\gamma_i>\alpha_i=\beta_i$ and $\gamma_j\geq min\{\alpha_j,\beta_j\}$ for each $j\neq i$, with equality holding if $\alpha_j\neq\beta_j$.

Given a semigroup $S$ of $\Z^r$ and a subset $E$ of $\Z^r$, with $S+E\subset E$ and such that $S$ and $E$ have the above properties $E_0$, $E_1$ and $E_2$, then we call $S$ a \emph{good semigroup} and $E$ a \emph{good semigroup ideal} of $S$. 

So, if $S=E(R)$ and $E=E(\I)$, where $\I$ is a regular fractional ideal, then $E$ a good semigroup ideal of $S$.

For a good semigroup ideal $E$, combining Properties $E_0$ and $E_1$, it follows that there exists a unique $m=m_E=\min(E)$.

On the other side, one has that if $\beta, \beta'\in E$ are such that  $\beta +\N^r \subset E$ and $\beta' +\N^r \subset E$, then
$\min(\beta,\beta') +\N^r \subset E$. This guarantees that there is a unique least element $\gamma\in E$ with the property that $\gamma+\N^r \subset E$. This element is called the \emph{conductor} of $E$ and denoted by $c(E)$. In particular, when $E=E(\mathcal I)$ for a fractional ideal $\I$, then we write $c(\I)$ for $c(E)$.


For a good semigroup ideal $E$, we will use the following notation:
\[\f(E) =c(E)-e, \quad \text{where} \ e=(1,\ldots,1),\]
which is called the \emph{Frobenius vector} of $E$. For $J\subset I$, we define $e_J$ the vector such that $pr_{\{i\}}(e_J)=1$ if $i\in J$ and $pr_{\{i\}}(e_J)=0$ if $i\not\in J$ and define $e_i=e_{\{i\}}$.  If $E=E(\I)$, we write $\f(\I)$ intead of $\f(E)$.

Since completion and value sets of fractional ideals are compatible (cf. \cite[\S 1]{Da} or \cite[Theorem 3.19]{KST}), we may assume that $R$ is complete with respect to the $\mathfrak M$-adic topology. In this case, the number $r$ of discrete valuation rings of $Q$ over $R$ coincides with the number of minimal primes of $R$. 


A fundamental notion in our context, which we define below, is that of a \emph{fiber} of an element $\alpha=(\alpha_1, \ldots, \alpha_r)\in E \subset \Z^r$ with respect to a subset $J=\{j_1<\cdots <j_s\}\subset I$. We define $pr_J(\alpha)=(\alpha_{j_1},\ldots,\alpha_{j_s})$.

Given $E\subset \Z^r$, $\alpha\in\mathbb{Z}^r$ and $\emptyset \neq J\subset I$, we define: 	
\[
\begin{array}{lcl}
F_J(E,\alpha)&=&\{\beta\in E;\; pr_J(\beta)=pr_J(\alpha) \ \text{and} \ \beta_i > \alpha_i, \forall i\in I\setminus J\}, \\
 \overline{F}_J(E,\alpha)&=&\{\beta\in E;\; pr_J(\beta)=pr_J(\alpha), \  \text{and} \ \beta_i \geq  \alpha_i, \forall i\in I\setminus J\}, \\ 
F(E,\alpha)&=& \bigcup_{i=1}^rF_i(E,\alpha), \quad \text{where} \ F_i(E,\alpha)=F_{\{i\}}(E,\alpha).
\end{array}
\] 
	
The last set, above, will be called the \emph{fiber} of $\alpha$.  Notice that $F_I(E,\alpha)=\{\alpha\}$, if and only if $\alpha\in E$, otherwise $F_I(E,\alpha)=\emptyset$.\smallskip

The importance of this notion may be seen, for example, by the following result (cf. \cite{Da}): Up to a multiplicative unit in $\widetilde{R}$,  there is a unique canonical ideal $\omega^0$ such that $R\subset \omega^0 \subset\widetilde{R}$ and $E(\omega^0)=E^0$, where 
\begin{equation}\label{valuecanonical0}
E^0=\{\alpha \in \Z^r; \ F(E(R),\f(R)-\alpha)=\emptyset\}.
\end{equation}
Notice that $\f(\omega^0)=\f(R)$ (cf. \cite[Lemma 5.10]{KST}).

Since any canonical ideal $\omega$ of $R$ is a multiple of $\omega^0$ by a unit $u$ in $Q$, then $E(\omega)$ is a translation of $E(\omega^0)$ by $v(u)=\f(\omega)-\f(R)$. This leeds to  the following:
\[ 
E(\omega)=\{\alpha \in \Z^r; \ F(E(R),\f(\omega)-\alpha)=\emptyset\}.
\]

A property related to the fibers of the frobenius of a good semigroup ideal $E$ is that (cf. \cite[Lemma 4.1.10  ]{KST}):
\begin{equation}\label{frobfiber}
F(E,\f(E))=\emptyset.
\end{equation}

We will use later the following remark that follows readily from the definitions of fibers.
\begin{rem} \label{fibrafechada}  If $E\subset \Z^r$, $\alpha\in \Z^r$, $J\subset I$ and $J^c=I\setminus J$, then one has
\[
F_J(E,\alpha)= \bigcap_{i\in J} \overline{F}_i(E,\alpha+e_{J^c})
\]
In particular, for $J=\{i\}$ one has that
\[
F_i(E,\alpha)=\overline{F}_i(E,\alpha+e-e_i).
\]\end{rem}

Another fundamental notion is that of maximal points of good semigroup ideals.

 Let $E\subset \Z^r$ and $\alpha\in E$. We will say that $\alpha$ is a \emph{maximal} point of $E$ if $F(E,\alpha)=\emptyset$. 
\smallskip

This means that there is no element in $E$ with one coordinate equal to the corresponding coordinate of $\alpha$ and the other ones bigger.	

When $E$ is a good semigroup ideal, since it has a minimum $m_E$ and a conductor $\gamma=c(E)$, one has immediately that all maximal elements of $E$ are in the limited region
\[
\{(x_1,\ldots,x_r)\in \Z^r; \ m_{E}\leq x_i < \gamma_i, \ \ i=1,\ldots,r\}.
\]

This implies that $E$ has finitely many maximal points.

Next, we will describe some special types of maximal points that may occur in a good semigroup ideal $E$. For $\alpha\in E$, let
$$\begin{array}{l} p(E,\alpha)=\max\{n; F_J(E,\alpha)=\emptyset,\forall J\subset I,\#J\leq n\}, \ \text{and}\\ \\
q(E,\alpha)=\min\{n; F_J(E,\alpha)\neq\emptyset,\forall J\subset I,\#J\geq n\}.
\end{array}$$

Notice that $p(E,\alpha)<q(E,\alpha)$, and that $\alpha\in E$ if and only if $q(E,\alpha) \leq r$. Also, $\alpha\in E$ is a maximal point of $E$ if, and only if, $p(E,\alpha) \geq 1$.

Let $\alpha$ be a maximal point of $E$. We will call $\alpha$ an \emph{absolute maximal}, if $F_J(E,\alpha)=\emptyset$ for every $J\subset I$, $J\neq I$; that is, if and only if $p(E,\alpha)=r-1$. We will call $\alpha$ a \emph{relative maximal}, if $F_J(E,\alpha)\neq\emptyset$, for every $J\subset I$ with $\#J\geq2$; that is, $p(E,\alpha)=1$ and $q(E,\alpha)=2$. 
If $p=p(E,\alpha)\geq 1$ and $q=q(E,\alpha)\leq r$, we call $\alpha$ a \emph{maximal of type} $(p,q)$.
\smallskip

Delgado in \cite[Theorem 1.5]{D87} showed that $E(R)$ is determined recursively, in a combinatorial sense, by its set of relative maximals.  With essentially the same proof, one may show the same for any good semigroup ideal $E$.

\section{Symmetry}

The central results in this section will be several generalizations of results in \cite{D88}, \cite{CDK}, \cite{Pol17}, \cite{Pol18} and \cite{KST}, which establish some symmetry among  $E(\J\colon \I)$ and $E(\I)$ mediated by $E(\J)$ and among their maximal points. 

\begin{lema}\label{fibra}
Let $\I$ and $\J$ be any fractional ideals of $R$, then one has 
\[
E(\J:\I) \subseteq \{\beta\in\mathbb{Z}^r; F(E(\I),\f(\J)-\beta)=\emptyset\}.
\]
\end{lema}
\pf Let $\beta\in E(\J:\I)$ and suppose that $F(E(\I),\f(\J)-\beta)\neq\emptyset$. Then there exist $i\in I$ and $\alpha\in E(\I)$ such that $\alpha_i=\f(\J)_i-\beta_i$ and $\alpha_j>\f(\J)_j-\beta_j$, for $j\neq i$, $j\in I$.  From this it follows that $\alpha_i+\beta_i=\f(\J)_i$ and $\alpha_j+\beta_j>\f(\J)_j$, for all $i\neq j$ and since from (\ref{$E+E^*$}) we know that $\alpha+\beta\in E(\J)$, we get that $\alpha+\beta\in F_i(E(\J),\f(\J))$, which is a contradiction, because from (\ref{frobfiber}) we know that $F(E(\J),\f(\J))=\emptyset$.
\cqd

\begin{teo}\label{J:I}
Let $\I$ and $\J$ be fractional regular ideals of $R$. The following are equivalent:
\begin{enumerate}[\rm i)]
	\item $\J$ is a canonical ideal;
	\item  $E(\J\colon\I)=\{\beta\in\Z^r; F(E(\I),\f(\J)-\beta)=\emptyset\}$, for all $\I$.
\end{enumerate}
\end{teo}
\pf  i) $\Rightarrow$ ii) Since $\J$ is a canonical ideal, from \cite[Proposition 5.18]{KST} we know that $E(\J)=\alpha_0+E(\omega^0)$, where $\alpha_0= \f(\J)-\f(\omega^0)\in\Z^r$. From \cite[Theorem 5.27]{KST}, for any fractional ideal $\I$, we have
$$\begin{array}{lcl}
E(\J\colon\I)&=&E(\J)-E(\I)\\
&=&(\alpha_0+E(\omega^{0}))-E(\I) \\
&=& \alpha_0+(E(\omega^{0})-E(\I))  \\
&=&\alpha_0+\{\beta\in\Z^r, F(E(\I),\f(\omega^{0})-\beta)=\emptyset\}\\
&=&\{\beta\in\Z^r;F(E(\I),\f(\J)-\beta)=\emptyset\},
\end{array}$$
where the fourth equality follows from \cite[Lemma 5.16]{KST}.\smallskip

\noindent ii) $\Rightarrow$ i) Suppose that $E(\J\colon\I)=\{\beta\in\Z^r; F(E(\I),\f(\J)-\beta)=\emptyset\}$, for all fractional ideal $\I$. In particular, for $\I=R$ we have that
\begin{equation} \label{E(J)}
E(\J)=E(\J\colon R)=\{\beta\in\Z^r; F(E(R),\f(\J)-\beta)=\emptyset\}.
\end{equation}

We will show that $E(\J)=\alpha+E(\omega^{0})$, for $\alpha= \f(\J)-\f(\omega^{0})$, which will imply, by \cite[Proposition 5.18 and Theorem 5.25.]{KST}, that $\J$ is a canonical ideal.

Let $\beta\in E(\omega^{0})$, from (\ref{valuecanonical0}), we get 
$$
F(E(R),\f(\J)-(\alpha+\beta))= F(E(R),\f(\omega^{0})-\beta)=\emptyset.
$$
Hence, from (\ref{E(J)}), $\alpha+\beta \in E(\J)$, which yields $\alpha + E(\omega^{0})\subset E(\J)$.

Let now $\beta\in E(\J)$ and write $\beta=\alpha+\beta'$, with $\beta'\in\Z^r$. Since $\beta\in E(\J)$,  we have 
$$
\emptyset= F(E(R),\f(\J)- \beta)=  F(E(R),\f(\omega^{0})-\beta'), 
$$
hence, from (\ref{valuecanonical0}), $\beta'\in E(\omega^{0})$; and therefore $E(\J)\subset \alpha+E(\omega^{0})$, concluding the proof of our result.\cqd

As mentioned in the proof of the above theorem, the implication (i) $\Rightarrow$ (ii) was proved in the particular case in which $\J=\omega^0$ in \cite[Lemma 5.16 and Theorem 5.27]{KST} (see also \cite[Theorem 2.15]{Pol18}), while the converse is new. 

This leeds to the following result:
 
\begin{cor} \label{gor} The following are equivalent:
\begin{enumerate}[\rm i)]  
\item $E(R\colon \I)= \{ \beta\in \Z^r; \, F(E(\I),\f(R)-\beta)=\emptyset\}$, for all fractional ideal $\I$;
\item $R$ is Gorenstein.
\end{enumerate}
\end{cor}
\pf If one has equality for all fractional ideal $\I$, then for $\I=R$ one has that
\[
E(R)=E(R\colon R)=\{ \beta\in \Z^r; \, F(E(R),\f(R)-\beta)=\emptyset\},
\]
which says that $E(R)$ is symmetric, hence $R$ is Gorenstein (cf. \cite[Proposition 5.29]{KST}). Conversely, if $R$ is Gorenstein, then $R=\omega^0$ is a canonical ideal, and the result follows from Theorem \ref{J:I}.
\cqd

For $\alpha\in \Z^r$ and $\I$ a fractional ideal of $R$, we define 
\[
\I(\alpha)=\{h\in \I^{reg}; \ v(h)\geq \alpha\}.
\]

We denote by $\ell(M)$ the length of an $R$-module $M$. We have the following result.

\begin{lema}[{\rm \cite[Proposition 2.2]{Da} or \cite[Lemma 3.18]{KST}}] \label{lema}
If $\alpha\in\mathbb{Z}^r$, then we have
$$\ell \left(\dfrac{\I(\alpha)}{\I(\alpha+e_i)}\right)=\left\{
\begin{array}{ll}
1, & \ if \  \overline{F}_i(E(\I),\alpha)\neq \emptyset, \\ \\
0, & \ \text{otherwise},
\end{array}
\right.$$
\end{lema}

The following theorem generalizes \cite[Theorem 3.6]{CDK}.

\begin{teo}\label{l_E}
	Let $\J$ and $\I$ be  fractional ideals of $R$ and let $\alpha,\beta\in\mathbb{Z}^r$, with $\alpha+\beta=c(\J)$. Then
	\begin{equation} \label{leq1}	
	\ell\left(\frac{\I(\alpha)}{\I(\alpha+e_i)}\right)+\ell\left(\frac{(\J\colon \I)(\beta-e_i)}{(\J\colon \I)(\beta)}\right)\leq1, \ \ \text{for every} \ i\in I ,\end{equation}
	with equality holding for every $\alpha,\beta$ such that $\alpha+\beta=c(\J)$ and for every  fractional ideal $\I$ if and only if $\J$ is a canonical ideal.
	\end{teo}
\pf Since by Lemma \ref{lema} each summand in (\ref{leq1}) is less than or equal to $1$, it is sufficient to show that they are not both equal to $1$.

Suppose by reductio ad absurdum that both summands in (\ref{leq1}) are equal to $1$. From Lemma \ref{lema}, it follows that  $\overline{F}_i(E(\I),\alpha)\neq\emptyset$ and $\overline{F}_i(E(\J:\I),\beta-e_i)\neq \emptyset$. 
	Take $\theta$ in the first of the above two sets and $\theta'$ in the second one, then according to (\ref{$E+E^*$}) we have $\theta+\theta'\in E(\J)$; even more, we have that $\theta+\theta'\in F_i(E(\J),\f(\J))$, because $\theta_i+\theta'_i=\f(\J)_i$ and $\theta_j+\theta'_j>\f(\J)_j$ for all $j\neq i$, which is a contradiction, since $F(E(\J),\f(\J))=\emptyset$.

Assuming that the equality holds in (\ref{leq1}), we will show that $\J$ is a canonical ideal.


Notice that, in view of Lemma \ref{lema}, equality in (\ref{leq1}) is equivalent to  
\begin{equation}\label{equiv}
 \overline{F}_i(E(\I),\alpha) = \emptyset \ \Iff \ \overline{F}_i(E(\J:\I),\beta-e_i) \neq \emptyset,  \ \forall i\in I.
\end{equation}

From the inclusion in Lemma \ref{fibra} we know that 
\begin{equation}\label{a}
E(\J\colon \I) \subset \{\beta\in \Z^r; \ F(E(\I),\f(\J)-\beta)=\emptyset\}.
\end{equation}

On the other hand, suppose that $\beta$ is such that $F(E(\I),\f(\J)-\beta)=\emptyset$, hence for all $i\in I$, $F_i(E(\I),\f(\J)-\beta)=\emptyset$. Then, 
from Remark \ref{fibrafechada}, it follows that $\overline{F}_i(E(\I),\f(\J)+e-\beta-e_i)=\emptyset$, for all $i\in I$. Now, from (\ref{equiv}), we get 
$\overline{F}_i(E(\J\colon \I),\beta)\neq \emptyset$, for all $i\in I$, which implies that $\beta\in E(\J\colon \I)$. 

Hence, we have shown that the inclusion in (\ref{a}) is an equality, therefore, from Theorem \ref{J:I}, we have that $\J$ is a canonical ideal.

Let us assume conversely that $\J$ is a canonical ideal. From Theorem \ref{J:I} we know that 
\[
\beta\notin E(\J\colon \I) \ \Leftrightarrow \ F(E(\I),\f(\J)-\beta)\neq \emptyset.
\]

To conclude the proof of this part of the theorem, it is clearly enough to show that 
\[
\forall \ i\in I, \ \ \overline{F}_i(E(\J\colon \I),\beta-e_i)=\emptyset \ \Longrightarrow \ \overline{F}_i(E(\I),\alpha)\neq \emptyset.
\]

Suppose that $\overline{F}_i(E(\J\colon \I),\beta-e_i)=\emptyset$, for some $i$. This implies that  $\beta-e_i\notin E(\J\colon \I)$, so, from Theorem \ref{J:I} we get $F(E(\I),\f(\J)-\beta +e_i)\neq \emptyset$. Take now $\theta \in F_i(E(\I),\f(\J)-\beta +e_i)$, for some $i\in I$. So, $\theta_i=\f(\J)_i-\beta_i +1=c(\J)_i-\beta_i=\alpha_i$ and $\theta_j\geq \f(\J)_j-\beta_j+1=c(\J)_j-\beta_j=\alpha_j$, for all $j\neq i$. So, we have that $\theta\in \overline{F}_i(E(\I),\alpha)$, hence this set is nonempty.
\cqd

We have the following consequence of the above theorem.

\begin{cor}\label{dimGo} For every fractional ideal $\I$ of $R$ and every $\alpha,\beta\in \Z^r$, with $\alpha+\beta=c(R)$, one has
\[
\ell\left(\frac{\I(\alpha)}{\I(\alpha+e_i)}\right)+\ell\left(\frac{(R\colon \I)(\beta-e_i)}{(R\colon \I)(\beta)}\right) =1 \ \text{for every} \ i\in I 
\]
if and only if $R$ is Gorenstein.
\end{cor}
	
Let $\I$ and $\J$ be fractional ideals of $R$ and let $\alpha\in \Z^r$. Let us define
\[
\rho_{\J}(\I,\alpha)=p(E(\I),\alpha)+q(E(\J\colon \I),\f(\J)-\alpha)-1.
\]

The following theorem generalizes \cite[Theorem 5.3]{CDK}. Recall that we defined $J^c$ as being the complement in $I$ of any subset $J$ of $I$.

\begin{teo}\label{p+q}
	For any fractional ideals $\I$ and $\J$ of $R$ and for any $\alpha\in\mathbb{Z}^r$, we have
\begin{equation}\label{p,q}
\rho_{\J}(\I,\alpha) \geq r.
\end{equation}
Moreover, equality holds in (\ref{p,q}), for every fractional ideal $\I$ of $R$ and every $\alpha\in\mathbb{Z}^r$ if, and only if, $\J$ is a canonical ideal.
\end{teo}
\pf
Suppose that $q(E(\J\colon\I),\f(\J)-\alpha)=r-n+1$, from the definition of $q$,  we know that
\begin{equation}\label{K}
F_{K}(E(\J\colon\I),\f(\J)-\alpha)\neq\emptyset, \ \forall K \subset I, \ \#K \geq r-n+1.
\end{equation}

If we take any $J\subset I$ with $\#J\leq n$ and let $K=J^c\cup\{i\}$, where $i\in J$ is fixed, then $\#K\geq r-n+1$. Hence from (\ref{K}) we get easily that 
$$\overline{F}_{K}(E(\J\colon\I),\f(\J)-\alpha+e_{K^c})\neq\emptyset.
$$

Since $e_{K^c}=e-e_{K}=e-(e_{J^c}+e_i)$, we get that  
$$ \overline{F}_{K}(E(\J\colon\I),c(\J)-\alpha-e_{J^c}-e_i)\neq\emptyset.$$

This implies that
$$\overline{F}_{i}(E(\J\colon\I),c(\J)-\alpha-e_{J^c}-e_i)\neq\emptyset, $$
that, from Theorem \ref{l_E} and Lemma \ref{lema}, implies that 
$$\overline{F}_{i}(E(\I),\alpha+e_{J^c})=\emptyset, \ \forall i\in J,$$
which in view of Remark \ref{fibrafechada}, implies that $F_J(E(\I),\alpha)=\emptyset$. 

So, we have shown, for any $J\subset I$, with $\#J\leq n$, that $F_J(E(\I),\alpha)=\emptyset$. Hence it follows that $p(E(\I),\alpha)\geq n$, and consequently, $\rho_{\J}(\I,\alpha)\geq r$. \smallskip

Now, suppose that $\J$ is a canonical ideal. To prove equality in (\ref{p,q}) holds, we must show that $p(E(\I),\alpha)=n$. Suppose by reductio ad absurdum that $p(E(\I),\alpha)\geq n+1$. Then, from the definition of $p$, we have that 
$$F_J(E(\I),\alpha)=\emptyset, \ \forall J\subset I, \ \text{with} \ \#J=n+1,
$$
 which implies that
	\begin{equation}\label{J^c}
	\overline{F}_i(E(\I),\alpha+e_{J^c})=\emptyset, \ \forall i\in J, \ \text{with}\ \#J\leq n+1,
	\end{equation}
	because, otherwise, we would have for some $i\in J$ that $\overline{F}_i(E(\I),\alpha+e_{J^c})\neq \emptyset$. Take $\theta$ is this last nonempty set, then $\theta_i=\alpha_i$, $\theta_j \geq \alpha_j$, $\forall j\in J$ and $\theta_l>\alpha_l$, $\forall l\not\in J$. Let $J'$ be the subset of elements $j\in J$ such that $\theta_j =\alpha_j$, hence $\theta \in F_{J'}(E(\I),\alpha)$, which implies that $F_{J'}(E(\I),\alpha)\neq  \emptyset$, with $\#J' \leq n+1$, contradicting the fact that $p(E(\I),\alpha)\geq n+1$.
	
For any $K\subset I$, with $\#K=r-n$, define the set $J=K^c\cup \{i\}$, where $i\in K$. Since $J$ has $n+1$ elements, it follows from (\ref{J^c}) that 
	$$\overline{F}_i(E(\I),\alpha+e_{J^c})=\emptyset, \ \forall i\in J.$$
	
	Since, $\J$ is a canonical ideal, from Theorem \ref{l_E} and Lemma \ref{lema}, it follows that
	$$\overline{F}_i(E(\J\colon\I),\f(\J)-\alpha+e_{K^c})\neq\emptyset,$$
	and since, $i$ was any element of $K$, we have that
	$$\overline{F}_i(E(\J\colon\I),\f(\J)-\alpha+e_{K^c})\neq\emptyset,\ \forall K\subset I, \ \#K=r-n, \forall i\in K.$$
	
	For every $i\in K$, take $\theta^i\in \overline{F}_i(E(\J\colon\I),\f(\J)-\alpha+e_{K^c})$, then $\theta^i_i=\f(\J)_i-\alpha_i$. $\theta^i_k\geq \f(\J)_k-\alpha_k$, for $k\in K$ and $\theta^i_j>\f(\J)_j-\alpha_j$, for $j\not\in K$. If we take $\theta=\min\{\theta^i; \ i\in K\}$, it follows that $\theta \in F_K(E(\J\colon\I),\f(\J)-\alpha)$, hence 
\[
F_K(E(\J\colon\I),\f(\J)-\alpha)\neq \emptyset, \ \forall K\subset I, \ \#K=r-n,
\]
which contradicts the fact that $q(E(\J\colon\I),\f(\J)-\alpha)=r-n+1$. Therefore, we must have $p(E(\I),\alpha)=n$.\smallskip

Now, assume that we have equality in (\ref{p,q}). Let $\I$ be a fractional ideal of $R$ and let $\alpha\in\mathbb{Z}^r$. If $\overline{F}_i(E(\I),\alpha)=\emptyset$, for some $i\in I$, then, from \cite[Lemma 4.7]{CDK}, there exists $\beta$ with $\beta_i=\alpha_i$ and $\beta_j<\alpha_j$ for every $j\neq i$, such that $F(E(\I),\beta)=\emptyset$. From this last condition, we get that $p(E(\I),\beta)\geq1$, so, from the equality in (\ref{p,q}), we get that  $q(E(\J\colon\I),\f(\J)-\beta)\leq r$, which means that $\f(\J)-\beta\in E(\J\colon\I)$. This implies that $\overline{F}_i(E(\J\colon\I),\f(\J)-\beta)\neq\emptyset$. Now,  since $\f(\J)_i-\beta_i=\f(\J)_i-\alpha_i$ and $\f(\J)_j-\beta_j> \f(\J)_j-\alpha_j$, it follows that $\emptyset \neq \overline{F}_i(E(\J\colon\I), \f(\J)-\beta) \subset \overline{F}_i(E(\J\colon\I), \f(\J)-\alpha)$, hence this last set is nonempty. So, we proved that equality holds in (\ref{leq1}), which, by Theorem \ref{l_E}, implies that $\J$ is a canonical ideal. \cqd

This leads immediately to the following result:

\begin{cor} The following two conditions are equivalent:

\noindent i) $\rho_R(\I,\alpha) = r$, for all fractional ideal $\I$ and all $\alpha\in \Z^r$;

\noindent ii) $R$ is Gorenstein.
\end{cor}

The following result will generalize \cite[Theorem 2.10]{D87}. 

\begin{teo}[\textbf{Symmetry of maximals}]\label{ed3} Let $\I$ and $\J$ be  fractional ideals of $R$. Suppose that $\alpha\in E(\I)$ and $\f(\J)-\alpha\in E(\J\colon \I)$. Then $\alpha$ is maximal of $E(\I)$ if  and only if $\f(\J)-\alpha$ is maximal of $E(\J\colon \I)$. Moreover, if $\alpha$ is a maximal of type $(p,q)=(p(E(\I),\alpha),q(E(\I),\alpha))$ then $\f(\J)-\alpha$ is maximal of type $(p',q')$,   where

$p'=\rho_\J(\J\colon(\J\colon \I),\f(\J)-\alpha)+1-q(E(\J\colon(\J\colon \I)),\alpha), \ \text{and}$

$q'=\rho_\J(\I,\alpha)+1-p(E(\I),\alpha).$
\end{teo}
\pf Suppose that $\alpha$ is a maximal of $E(\I)$, then $p(E(\I),\alpha)\geq1$. From Theorem \ref{p+q} we have
\[\rho_{\J}(\J\colon\I,\f(\J)-\alpha)=p(E(\J\colon\I),\f(\J)-\alpha)+q(E(\J\colon(\J\colon\I)),\alpha)-1\geq r.\]

Since $\I\subset \J\colon(\J\colon\I)$, by the definition of the number $q$ we have
\[q(E(\J\colon(\J\colon\I)),\alpha)\leq q(E(\I),\alpha), \ \text{for any} \ \alpha\in\Z^r.\]

Hence, 
$$
\begin{array}{lcr}
r&\leq &p(E(\J\colon\I),\f(\J)-\alpha)+q(E(\J\colon(\J\colon\I)),\alpha)-1 \\
  &\leq & p(E(\J\colon\I),\f(\J)-\alpha)+q(E(\I),\alpha)-1,
	\end{array}$$
so $p(E(\J\colon\I),\f(\J)-\alpha)\geq1$ and, since $\f(\J)-\alpha\in E(\J\colon \I)$, it follows that $\f(\J)-\alpha$ is a maximal of $E(\J\colon \I)$.

The proof of the converse of this statement is completely analogous. \smallskip

Furthermore, since $p'=p(E(\J\colon\I),\f(\J)-\alpha)$, we have 
\[ p'=p(E(\J\colon\I),\f(\J)-\alpha)=\rho_\J(\J\colon(\J\colon \I),\f(\J)-\alpha)+1-q(E(\J\colon(\J\colon \I)),\alpha);
\]
and since, for $\alpha\in E(\I)$, we have
\[
\rho_{\J}(\I,\alpha)=p(E(\I),\alpha)+q(E(\J\colon\I),\f(\J)-\alpha)-1,
\]
then 
\[q'=q(E(\J\colon\I),\f(\J)-\alpha) =\rho_{\J}(\I,\alpha)+1-p(E(\I),\alpha).
\]
\cqd

If $\J$ is a canonical ideal, do not have to assume that both $\alpha\in E(\I)$ and $\f(\J)-\alpha\in E(\J\colon \I)$, since in this  case, $$\alpha \ \text{is maximal of } \ E(\I) \ \Iff \ \f(\J)-\alpha \ \text{is maximal of} \ E(\J\colon \I).$$ Also, $\J\colon(\J\colon \I)=\I$ and $\rho_{\J}(\I,\alpha)=r$. Taking this into account, we get the following result:

\begin{cor} Suppose that $\J$ is a canonical ideal, then one has that $\alpha$ is maximal of $E(\I)$ if, and only if, $\f(\J)-\alpha$ is maximal of $E(\J\colon \I)$. Moreover, $\alpha$ is a maximal of type $(p,q)$ if and only if $\f(\J)-\alpha$ is a maximal of type $(r+1-q,r+1-p)$.
\end{cor}


\begin{thebibliography}{X}
\bibitem{Ap}\textsc{Ap\'ery, R.}; \textit{Sur les branches superlin\'eaires des courbes alg\'ebriques}, C.R.A.S. Paris, vol. 222, pp 1198-1200 (1946).

\bibitem{CDK} \textsc{Campillo, A.}; \textsc{Delgado de la Mata, F.}; \textsc{Kiyek, K.}, \textit{Gorenstein property and symmetry for one-dimensional local Cohen-Macaulay rings.}, Manuscripta Math. 83, pp 405-423 (1994).

\bibitem{Da} \textsc{D'anna, M.}, \textit{The Canonical Module of a One-dimensional Reduced Local Ring}, Communications in Algebra, 25, pp 2939-2965 (1997).

\bibitem{D87} \textsc{Delgado de la Mata, F.}, \textit{The semigroup of values of a curve singularity with several branches}, Manuscripta Math. 59, pp 347-374 (1987).

\bibitem{D88} \textsc{Delgado de la Mata, F.}, \textit{Gorenstein curves and symmetry of the semigroup of values}, Manuscripta Math. 61, pp 285-296 (1988).

\bibitem{Ga} \textsc{Garcia, A.}, \textit{Semigroups associated to singular points of plane curves},
J. Reine. Angew. Math. 336, 165-184 (1982).



\bibitem{ja} \textsc{J\"ager, J.}, \textit{L\"angeberechnungen und kanonische Ideale in eindimensionalen Ringen}, Arch. Math. 29, pp 504-512 (1977).

\bibitem{KST} \textsc{Korell, P.; Schulze, M.; Tozzo, L.}, \textit{Duality on value semigroups}, J. of Comm. Alg. V. 11, no. 1, pp 81-129 (2019).  


\bibitem{Ku} \textsc{Kunz, E.}, \textit{The value semigroup of a one-dimensional Gorenstein ring}, Proc. Amer. Math. Soc, 25, pp 748-751 (1970).


\bibitem{Pol17} \textsc{Pol, D.}, \textit{On the values of logarithmic residues along curves},  Annales de l'Institut Fourier, 68 no. 2, pp 725-766 (2018).

\bibitem{Pol18} \textsc{Pol, D.}, \textit{Symmetry of maximals for fractional ideals of curves}, To appear on J. of Comm. Alg. (arXiv:1802.07901).

\bibitem{W} \textsc{Waldi, R.}, \textit{Wertehalbgruppe und Singularit\"at einer ebenen algebraischen Kurve}, Dissertation. Regensburg (1972).
\bibitem{Za1} \textsc{Zariski, O.}, \textit{Studies in Equisingularity I}, Amer. Jour. Math. 87, pp 507-536 (1965).


\end{thebibliography}
\end{document}